\documentclass{emsprocart1}
\usepackage{mathrsfs,amsmath}
\usepackage[all]{xy}

\newcommand{\remark}[1]{\subsection{Remark} {#1}}
\newcommand{\theorem}[2]{\subsection{Theorem #1} \itshape{#2}}
\newcommand{\lemma}[2]{\subsection{Lemma #1} \itshape{#2}}

\newcommand{\Z}{\mathbb{Z}}
\newcommand{\N}{\mathbb{N}}
\newcommand{\C}{\mathbb{C}}
\newcommand{\Q}{\mathbb{Q}}
\newcommand{\R}{\mathbb{R}}
\newcommand{\A}{\mathscr{A}}
\newcommand{\B}{\mathscr{B}}
\newcommand{\K}{\mathscr{K}}
\newcommand{\T}{\mathcal{T}}
\newcommand{\ootimes}{\widehat{\otimes}}
\newcommand{\G}{{\displaystyle\widetilde{G}}}
\title[TWISTED $K$-THEORY OLD AND NEW]{TWISTED $K$-THEORY\\OLD AND NEW}
\author[Max KAROUBI]{by Max KAROUBI}
\contact[max.karoubi@gmail.com]
{Max Karoubi\\
Universit\'e Paris 7\\
2 place Jussieu\\
75251 PARIS (FRANCE)}
\begin{document}
\begin{abstract}
Twisted $K$-theory has its origins in the author's PhD thesis \cite{K1} and in a
paper with P. Donovan \cite{DK}. The objective of this paper is to revisit the
subject in the light of new developments inspired by Mathematical Physics. See for
instance E. Witten \cite{Wi}, J. Rosenberg \cite{R}, C. Laurent-Gentoux, J.-L. Tu,
P. Xu \cite{TX} and M.F. Atiyah, G. Segal \cite{AS1}, among many authors. We also
prove some new results in the subject: a Thom isomorphism, explicit computations in
the equivariant case and new cohomology operations.
\end{abstract}
\begin{classification}
55N15, 19K99, 19L47, 46L80
\end{classification}

\begin{keywords}
$K$-theory, Brauer group, Fredholm operators, Poincar\'e duality, Clifford algebras,
Adams operations.
\end{keywords}
\section*{Some history and motivation about this paper}
The subject ``$K$-theory with local coefficients", now called ``twisted $K$-theory",
was introduced by P. Donovan and the author in \cite{DK} almost forty years
ago.\footnote{See the appendix for a short history of the subject.} It associates to
a compact space $X$ and a ``local coefficient system"
$$\alpha\in GBr(X)=\Z/2\times H^1(X;\Z/2)\times Tors(H^3(X;\Z))$$
an abelian group $K^\alpha(X)$ which generalizes the usual
Grothendieck-Atiyah-Hirzebruch $K$-theory of $X$ when we restrict  $\alpha$ being in
$\Z/2$ (cf. \cite{AH2}). This ``graded Brauer group" $GBr(X)$ has the following
group structure: if $\alpha = (\varepsilon, w_1, W_3)$ and $\alpha' =
(\varepsilon',w_1',W_3')$ are two elements, one defines the sum $\alpha + \alpha'$
as $(\varepsilon+ \varepsilon', w_1 +w_1',W_3 +W_3' + \beta(w_1\cdot w_1'))$, where
$\beta:H^2(X ; \Z/2)\to H^3(X; \Z)$ is the Bockstein homomorphism. With this
definition, one has a generalized cup-product\footnote{Strictly speaking, this
product is defined up to non canonical isomorphism ; see 2.1 for more details.}.
$$K^\alpha(X) \times K^{\alpha'}(X)\to K^{\alpha+\alpha'}(X)$$
The motivation for this definition is to give in $K$-theory a satisfactory Thom
isomorphism and Poincar\'e duality pairing which are analogous to the usual ones in
cohomology with local coefficients. More precisely, as proved in \cite{K1}, if $V$
is a real vector bundle on a compact space $X$ with a positive metric, the
$K$-theory of the Thom space of $V$ is isomorphic to a certain ``algebraic'' group
$K^{C(V)}(X)$ associated to the Clifford bundle $C(V)$, viewed as a bundle of
$\Z/2-$graded algebras. A careful analysis of this group shows that it depends only
on the class of $C(V)$ in $GBr(X)$, the three invariants being respectively the rank
of $V$ mod. 2, $w_1(V)$ and $\beta(w_2(V)) = W_3(V)$, where $w_1(V)$, $w_2(V)$ are
the first two Stiefel-Whitney classes of $V$. In particular, if $V$ is a
${}^c$spinorial bundle of even rank, one recovers a well-known theorem of Atiyah and
Hirzebruch. On the other hand, if $X$ is a compact manifold, it is well-known that
such a Thom isomorphism theorem induces a pairing between $K$-groups
$$K^\alpha(X)\times K^{\alpha'}(X)\to\Z$$
if $\alpha + {\alpha'}$ is the class of $-C(V)$ in $GBr(X)$, where $V$ is the tangent
bundle of $X$.\\\\
The necessity to revisit these ideas comes from a new interest in the subject
because of its relations with Physics \cite{Wi}, as shown by the number of recent
publications. However, for these applications, the first definition recalled above
is not complete since the coefficient system is restricted to the torsion elements
of $H^3(X; \Z)$. As it was pointed out by J. Rosenberg \cite{R} and later on by C.
Laurent, J.-L. Tu, P. Xu \cite{LTX}, M. F. Atiyah and G. Segal \cite{AS1}, this
restriction is in fact not necessary. In order to avoid it, one may use for instance
the Atiyah-J\"anich theorem \cite{A1}\cite{J} about the representability of
$K$-theory by the space of Fredholm operators
(already quoted in \cite{DK} for the cup-product which cannot be defined otherwise).\\
In the present paper, we would like to make a synthesis between different viewpoints
on the subject:  \cite{DK}, \cite{R}, \cite{LTX} and \cite{AS1} (partially of
course)\footnote{As it was pointed out to me by J. Rosenberg, one should also add
the following reference, in the spirit of \cite{DD}: Ellen Maycock Parker, The
Brauer group of graded continuous trace C*-algebras, Trans. Amer. Math. Soc. 308
(1988), Nr. 1, 115--132.}. We hope
have been ``pedagogical" in some sense to the non experts.\\\\
However, this paper is not just historical. It presents the theory with another
point of view and contains some new results. We extend the Thom isomorphism to this
more general setting (see also \cite{Ca}), which is important in order to relate the
``ungraded" and ``graded" twisted $K$-theories. We compute many interesting
equivariant twisted $K$-groups, complementing the basic papers \cite{LTX},
\cite{AS1} and \cite{AS2}. For this purpose, we use the ``Chern character" for
finite group actions, as defined by Baum, Connes, Kuhn and Slominska
\cite{BC}\cite{Kuh}\cite{SL}, together with our generalized Thom isomorphism. These
last computations are related to some previous ones \cite{K4} and to the work of
many authors. Finally, we introduce new cohomology operations which are
complementary to those defined in
\cite{DK} and \cite{AS2}.\\\\
We don't pretend to be exhaustive in a subject which has already many ramifications.
In an appendix to this paper, we try to give a short historical survey and a list of
interesting contributions of many authors related to the results quoted here.

\section*{General plan of the paper}

Let us first recall the point of view developed in \cite{DK}, in order to describe
the background material. We consider  a locally trivial bundle of $\Z/2$-graded central
simple complex algebras $A$, i.e. modelled on $M_{2n}(\C)$ or $M_{n}(\C)\times
M_{n}(\C)$, with the obvious gradings \footnote{up to graded Morita equivalence.}.
Then $\A$ has a well-defined class $\alpha$ in the group $GBr(X)$ (as introduced above). On
the other hand, one may consider the category of ``$\A$-bundles" , whose objects are
vector bundles provided with an $\A$-module structure (fibrewise). We call this
category $E^\A(X)$; the graded objects of this category are vector bundles which are
modules over $\A\ootimes C^{0,1}$, where $C^{0,1}$ is the Clifford algebra $\C\times
\C = \C[x]/(x^2 -1)$. The group $K^\alpha(X)$ is now defined as the ``Grothendieck
group'' of the forgetful functor
$$E^{\A\ootimes C^{0,1}}(X)\to E^\A(X)$$
We refer the reader to \cite{K1} p. 191 for this definition which generalizes the
usual Grothendieck group of a category. For our purpose, we make it quite explicit
at
the end of \S 1, using the concept of ``grading".\\\\
Despite its algebraic simplicity, this definition of  $K^\alpha(X)$ is not quite
satisfactory for various reasons. For instance, it is not clear how to define in a
simple way a ``cup-product"
$$K^\alpha(X)\times K^{\alpha'}(X)\to K^{\alpha+{\alpha'}}(X)$$
as mentioned earlier (even if $\alpha$ and  $\alpha'$ are in the much smaller group
$\Z/2$).\\\\
To correct this defect, a second definition may be given in terms of Fredholm operators
in a Hilbert space. More precisely, we consider graded Hilbert bundles $E$ which are also
graded $\A$-modules in an obvious sense, together with a continuous family of Fredholm
operators
$$D: E\to E$$
with the following properties:
\begin{enumerate}
 \item $D$ is self-adjoint of degree one
\item $D$ commutes with the action of $\A$ (in the graded sense)\\
\end{enumerate}
One gets an abelian semi-group from the homotopy classes of such pairs $(E, D)$,
with the addition rule
$$(E, D) + (E', D') = (E\oplus E', D\oplus D')$$
The associated group gives the second definition of $K^\alpha(X)$ which is
equivalent to the first one (see \cite{DK} p. 18 and \cite{K2} p. 88). We use also
the notation $K^\A(X)$ instead of $K^\alpha(X)$ where $\alpha$ is
the class of $\A$ in $GBr(X)$ when we want to be more explicit.\\\\
With this new viewpoint, the cup-product alluded to above becomes obvious. It is
defined by the following formula\footnote{See the appendix about the origin of such
a formula.}:

$$(E, D)\smile(E', D') = (E\ootimes E', D\ootimes 1 + 1\ootimes D')$$
where the symbol $\ootimes$ denotes the graded tensor product of bundles or
morphisms. It is a map from $K^\A(X)\times K^{\A'}(X)$ to $K^{\A\ootimes\A'}(X)$ and
therefore it induces a (non canonical) map from
$K^\alpha(X) \times K^{\alpha'}(X)$ to $K^{\alpha+\alpha'}(X)$.\\\\
For simplicity's sake, we have only considered \underline{complex} $K$-theory. We
could as well study the real case: one has to replace $GBr(X)$ by

$$GBrO(X) = \Z/8\times H^1(X; \Z/2)\times H^2(X; \Z/2)$$
If we take for the coefficient system $\alpha = n$ to be in $\Z/8$, we get the usual
groups $KO^n(X)$ as defined using Clifford algebras in \cite{K1} and \cite{K2} p. 88
(these groups
being written $\overline{K}^n$ in the later reference).\\\\
In this paper, we essentially follow the same pattern, but with bundles of infinite
dimension in the spirit of \cite{R}. As a matter of fact, all the technical tools
are already present in \cite{DK}, \cite{K2} and \cite{R}, for instance the Fredholm
operator machinery which is necessary to define the cup product. However, this paper
is not a rewriting of these papers, since we take a more synthetic view point and
have other applications in mind. For instance, the $K$-theory of Banach algebras
$K_n(A)$ and its graded version, denoted here by $GrK_n(A)$, are more systematically
used. On the other hand, since the equivariant $K$-theory and its relation with
cohomology have been studied carefully in \cite{CKRW}, \cite{LTX} and \cite{AS2}, we
limit ourselves to the applications in this case. One of them is the definition of
operations in twisted $K$-theory in the graded and ungraded
situations.\\\\
Here are the contents of the paper:

\begin{enumerate}

\item $K$-theory of $\Z/2$-graded Banach algebras
\item Ungraded twisted $K$-theory in the finite and infinite-dimensional cases
\item Graded twisted $K$-theory in the finite and infinite-dimensional cases
\item The Thom isomorphism
\item General equivariant $K$-theory
\item Some computations in the equivariant case
\item Operations in twisted $K$-theory.
\item[Appendix]. A short historical survey of twisted $K$-theory.
\end{enumerate}
\renewcommand\keywordname{Acknowledgments}

\begin{keywords} I would like to thank A. Adem, A. Carey, P. Hu, I. Kriz, C. Leruste,
V. Mathai, J. Rosenberg, J.-L. Tu, P. Xu and the referee for their remarks and
suggestions after various drafts of this paper.
\end{keywords}

\section{$K$-theory of $\Z/2$-graded Banach algebras.}

\subsection{}

Higher $K$-theory of real or complex Banach algebras $A$ is well-known (cf.
\cite{K1}
 or \cite{B} for instance). Starting from the usual Grothendieck group $K(A) = K_0(A)$,
 there are many equivalent ways to define ``derived functors'' $K_n(A)$, for $n\in\Z$, such
 that any exact sequence of Banach algebras

$$\xymatrix{0\ar[r]&A'\ar[r]&A\ar[r]&A''\ar[r]&0}$$
induces an exact sequence of abelian groups

$$\xymatrix@=.6cm{\dots\ar[r]&
K_{n+1}(A)\ar[r]& K_{n+1}(A")\ar[r]& K_n(A')\ar[r]& K_n(A)\ar[r]& K_n(A")\ar[r]&
\dots}$$
Moreover, by Bott periodicity, these groups are periodic of period 2 in the
complex case and 8 in the real case.

\subsection{}

The $K$-theory of $\Z/2$-graded Banach algebras $A$ (in the real or complex case) is
less well-known\footnote{This is of course included in the general $KK$-theory of
Kasparov which was introduced later than our basic references, at least for
C*-algebras.} and for the purpose of this paper we shall recall its definition which
is already present but not systematically used in \cite{K1} and \cite{DK}. We first
introduce $C^{p,q}$ as the Clifford algebra of $\R^{p+q}$ with the quadratic form
$$-(x_1)^2 -\dots -(x_p)^2 + (x_{p+1})^2 + \dots+ (x_{p+q})^2$$
It is naturally $\Z/2$-graded. If $A$ is an arbitrary $\Z/2$-graded Banach algebra,
$A^{p,q}$ is the graded tensor product $A\ootimes C^{p,q}$. For $A$ unital, we now
define the graded $K$-theory of $A$ (denoted by $GrK(A)$) as the $K$-theory of the
forgetful functor:
$$\phi:\mathcal{P}(A^{0,1})\to \mathcal{P}(A)$$
(see \cite{K1} or 1.4 for a concrete definition). Here $\mathcal{P}(A)$ denotes in
general the category of finitely generated projective left $A$-modules. One may
remark that the objects of $\mathcal{P}(A^{0,1})$ are graded objects of the category
$\mathcal{P}(A)$. The functor $\phi$ simply ``forgets" the grading. One should also
notice that $GrK(A^{p,q})$ is naturally isomorphic to $GrK(A^{p+1, q+1})$, since
$A^{p,q}$ is Morita equivalent to $A^{p+1, q+1}$ (in the graded
sense).\\
If $A$ is not unital, we define $GrK(A)$ by the usual method, as the kernel of the
augmentation map
$$GrK(A^+)\to GrK(k)$$
where $A^+$ is the $k$-algebra $A$ with a unit added ($k = \R$ or $\C$ , according
to the theory, with the trivial grading).

\subsection{}

There is a ``suspension functor" on the category of graded algebras, associating to
$A$ the graded tensor product $ A^{0,1}=A\ootimes C^{0,1}$. One of the fundamental
results\footnote{Strictly speaking, one has to replace the category
$\mathcal{C}^{p,q}$ with an arbitrary graded category. However, the proof of Theorem
2.2.2 in \cite{K1} easily extends to this case.} in \cite{K1} p. 210 is the fact
that this suspension functor is compatible with the usual one: in other words, we
have a well-defined isomorphism
$$GrK(A^{0,1})\equiv GrK(A(\R))$$
where $A(T)$ denotes in general the algebra of continuous maps $f(t)$ on the locally
compact space $T$ with values in $A$, such that $f(t)\to 0$ when $t$ goes to infinity. As
a consequence, we have an isomorphism between the following groups (for $n\geq 0$):
$$GrK(A^{0,n})\equiv GrK(A(\R^n))$$
which we call $GrK_n(A)$. More generally, we put $GrK_n(A) = GrK(A^{p,q})$ for $q -
p = n\in\Z$. These groups $GrK_n(A)$ satisfy the same exactness property as the
groups $K_n(A)$ above, from which they are naturally derived. They are of course
linked with them by the following exact sequence (for all $n\in\Z$):
$$\xymatrix{K_{n+1}(A\ootimes C^{0,1})\ar[r]&
K_{n+1}(A)\ar[r]& GrK_n(A)\ar[r]& K_n(A\ootimes C^{0,1})\ar[r]& K_n(A)}$$ In
particular, if we start with an ungraded Banach algebra $A$, we see that Bott
periodicity follows from these previous considerations, thanks to the periodicity of
Clifford algebras up to graded Morita equivalence: this was the main theme developed
in \cite{ABS}, \cite{W} and \cite{K1}, in order to give a more conceptual proof of
the periodicity theorems.

\subsection{}

When $A$ is unital, it is technically important to describe the group $GrK(A)$ in a
more concrete way. If $E$ is an object of $\mathcal{P}(A)$, a \underline{grading} of
$E$ is given by an involution $\varepsilon$ which commutes (resp. anticommutes) with
the action of the elements of degree 0 (resp. 1) in $A$. In this way, $E$ with a
grading $\varepsilon$ may be viewed as a module over the algebra $A\ootimes
C^{0,1}$. We now consider triples $(E, \varepsilon_1, \varepsilon_2)$, where
$\varepsilon_1$ and $\varepsilon_2$ are two gradings of $E$. The homotopy classes of
such triples obviously form a semi-group. Its quotient by the semi-group of
``elementary" triples (i.e. such that $\varepsilon_1 = \varepsilon_2$) is isomorphic
to $GrK(A)$.

\section{Ungraded twisted $K$-theory in the finite and infinite-dimensional cases}

\subsection{}

Let $X$ be a compact space and let us consider bundles of algebras $\A$ with fiber
$M_{n}(\C)$. As was shown by Serre \cite{G}, such bundles are classified by \v{C}ech
cocycles (up to \v{C}ech coboundaries):
$$g_{ji}: U_i\cap U_j\to PU(n)$$
where $PU(n)$ is $U(n)/ S^1$, the projective unitary group. In other words, the
bundle of algebras $\A$ may be obtained by gluing together the bundles of
C*-algebras $(U_i\times M_{n}(\C))$, using the transition functions $g_{ji}$ [Note
that $U(n)$ acts on $M_{n}(\C)$ by inner automorphisms and therefore induces an
action of $PU(n)$ on the algebra $M_{n}(\C)$]. The Brauer group of $X$ denoted by
$Br(X)$ is the quotient of this semi-group of bundles (via the tensor product) by
the following equivalence relation: $\A$ is equivalent to $\A'$ iff there exist
vector bundles $V$ and $V'$ such that the bundles of algebras $\A\ootimes END(V)$
and $A'\ootimes END(V')$ are isomorphic. It was proved by Serre \cite{G} that
$Br(X)$ is naturally isomorphic to
the torsion subgroup of $H^3(X;\Z)$.\\\\
The Serre-Swan theorem (cf. \cite{K1} for instance) may be easily translated in this
situation to show that the category of finitely generated projective $\A$-module bundles $E$
(as in \cite{DK}) is equivalent to the category $\mathcal{P}(A)$ of finitely
generated projective modules over $A = \Gamma(X,\A)$, the algebra of continuous
sections of the bundle $\A$. The key observation for the proof is that $E$ is a
direct factor of a ``trivial" $\A$-bundle; this is easily seen with finite
partitions of unity, since $X$ is compact. One should notice that if $\A$ is
equivalent to $\A'$ the associated categories $\mathcal{P}(A)$ and $\mathcal{P}(A')$
are equivalent. Note however that this equivalence is non canonical since $\A\otimes
END(V)$ and $\A'\otimes END(V')$ are not canonically isomorphic.

\begin{definition} The ungraded\footnote{We use the notation $K^{(\A)}(X)$, not to
be confused with the graded twisted $K$-group $K^\A(X)$ which will be defined in the
next section.} twisted $K$-theory $K^{(\A)}(X)$ is by definition the $K$-theory of
the ring $A$ (which is the same as the $K$-theory of the category $E^\A(X)$
mentioned in the introduction). By abuse of notation, we shall simply call it
$K{(\A)}$. We also define $K_n^{(\A)}(X)$ as the $K_n$-group of the Banach algebra
$\Gamma(X, \A)$. It only depends on the class of $\A$ in $Br(X)=Tors(H^3(X;\Z))$.
\end{definition}

\subsection{}
The key observation made by J. Rosenberg \cite{R} is the following: we can
``stabilize" the situation (in the C*-algebra sense), i.e. embed $M_{n}(\C)$ into
the algebra of compact operators $\K$ in a separable Hilbert space $H$, thanks to
the split inclusion of $\C^n$ in $l^2(\N)$. Now, a bigger group $PU(H) = U(H)/S^1$
is acting on $\K$ by inner automorphisms. If we take a \v{C}ech cocycle
$$g_{ji}:U_i\cap U_j\to PU(H)$$
we may use it to construct a bundle $\A$ of (non unital) C*-algebras with $\hbox{fiber
}\K$.\\\\
Let us now consider the commutative diagram
$$\xymatrix{
S^1\ar[r]\ar[d]& S^1\ar[d]\\
U(n)\ar[r]\ar[d]&U(H)\ar[d]\\
PU(n)\ar[r]&PU(H) }$$ Thanks to Kuiper's theorem \cite{Ku}, we remark that the
classifying space of $U(H)$ is contractible. Therefore, the classifying space
$BPU(H)$ of the topological group $PU(H)$, is a nice model of the Eilenberg-Mac Lane
space $K(\Z, 3)$ (compare with the well-known paper of Dixmier and Douady
\cite{DD}). Moreover, if we start with a finite-dimensional algebra bundle $\A$ over
$X$ with fiber $M_{n}(\C)$, the diagram above shows how to associate to $\A$ another
bundle of algebras $\A'$ with fiber $\K$, together with a C*-inclusion from $\A$ to
$\A'$. We note that the invariant $W_3(A )$ in $Br(X) = Tors(H^3(X; \Z))$ defined in
\cite{G} is simply induced by the classifying map from $X$ to $BPU(H)$ (which
factors through $BPU(n)$). In this finite example, it is an $n$-torsion class since
one has another commutative diagram

$$\xymatrix{
\mu_{n}\ar[d]\ar[r]&S^1\ar[d]\\
SU(n)\ar[d]\ar[r]&U(n)\ar[d]\\
PU(n)\ar@{=}[r]&PU(n)
}$$

\begin{theorem}{} The inclusion from $\A$ to $\A'$ induces an isomorphism

$$\xymatrix{K_r (\A) = K_r (\Gamma(X, \A))\ar[r]&K_r (\Gamma(X, \A')) = K_r (\A')}$$
where the $K_r$ define the classical topological $K$-theory of C*-algebras.
\end{theorem}

\begin{proof} The proof is classical for a trivial algebra bundle, since $\C$ is Morita equivalent to $\K$ (in the C*-algebra sense). It extends to the general case by a no less classical Mayer-Vietoris argument.
\end{proof}

{\begin{definition} Let now $\A$ be an algebra bundle with fiber $\K$ on a compact
space $X$ \underline{with structural group $PU(H)$}. We define $K^{(\A)}(X)$ (also
denoted by $K{(\A)}$) as the $K$-theory of the (non unital) Banach algebra
$\Gamma(X, \A)$. This $K$-theory only depends of the class $\alpha$ of $\A$ in
$H^3(X;\Z)$ and we shall also call it $K^{(\alpha)}(X)$ [Due to 2.4, this is a
generalization of definition 2.2].
\end{definition}}

\subsection{} Before treating the graded case in the next section, we would like to give
an equivalent definition of $K(\A)$ in terms of Fredholm operators, as was done in
\cite{DK} for the torsion elements in $H^3(X; \Z)$ and in \cite{AS1} for the general
case. The basic idea is to remark that $PU(H)$ acts not only on the C*-algebra of
compact operators in $H$, but also on the ring of bounded operators $End(H)$ and on
the Calkin algebra $End(H)/\K$ (with the norm topology\footnote{other topologies may
be also considered, see \cite{AS1}.})). Let us call $\B$ the algebra bundle with
fiber $End(H)$ associated to the cocycle defined in 2.1 and $\B/\A$ the quotient
algebra bundle. Therefore, we have an exact sequence of C*-algebras bundles

$$\xymatrix{0\ar[r]&\A\ar[r]&\B\ar[r]&\B/\A\ar[r]&0}$$
which induces an exact sequence for the associated rings of sections (thanks to a
partition of unity again)
$$\xymatrix{0\ar[r]&\Gamma(X, \A)\ar[r]&\Gamma(X,\B)\ar[r]&\Gamma(X, \B/\A))\ar[r]&0}$$
If $\B$ is trivial, it is well-known that the algebra of continuous maps from $X$ to
$End(H)$ has trivial $K_n$-groups because this algebra is flabby\footnote{A unital
Banach algebra $\Lambda$ is called flabby if there exists a continuous functor
$\tau$ from $\mathcal{P}(\Lambda)$ to itself such that $\tau+Id$ is isomorphic to
$\tau$. For instance, $\Lambda = End(H)$ is flabby since $\mathcal{P}(A)$ is
equivalent to the category of Hilbert spaces which are isomorphic to direct factors
in $H$; $\tau$ is then defined by the infinite Hilbert sum
$\tau(E)=E\oplus\dots\oplus E\oplus\dots$}. By a Mayer-Vietoris argument, it follows
that $K_n(\Gamma(X, \B))$ is also trivial. Therefore the connecting homomorphism
$$\xymatrix{K_1(\B/\A) = K_1(\Gamma(X, \B/\A))\ar[r]&K_0(\Gamma(X, \A))\equiv K^{(\A)}(X) = K(\A)}$$
is an isomorphism, a well-known observation in index theory.

\subsection{}
Let us now consider the elements of $\B$ which map onto $(\B/\A)^*$ via the map
$\pi$. These elements form a bundle of Fredholm operators on $H$ (the twist comes
from the action of $PU(H)$). This subbundle of $\B$ will be denoted by $Fredh(H)$.
Therefore, we have a principal fibration
$$\xymatrix{\Gamma(X, \A)\ar[r]&\Gamma(X, Fredh(H))\ar[r]^\pi&\Gamma(X, (\B/\A)^*)}$$
with contractible fiber the Banach space $\Gamma(X, \A)$ (this fibration admits a
local section thanks to Michael's theorem \cite{Mi}). Therefore, the space of
sections of $Fredh(H)$ has the homotopy type of $\Gamma(X, (\B/\A)^*)$. In
particular, the path components are in bijective correspondence via the map $\pi$.
The following theorem is a generalization of a well-known theorem of Atiyah and
J\"anich \cite{A1} \cite{J}:

{\theorem{{\cite{AS1}}} {The set of homotopy classes of continuous sections of the
fibration

$$\xymatrix{Fredh (H)\ar[r] &X}$$
is naturally isomorphic to $K^{(\A)}(X)$.}}

\begin{proof} As we have seen above, the two spaces $\Gamma(X, Fredh (H))$ and $\Gamma(X, (\B/\A)^*)$ have the same homotopy type. On the other hand, it is a well-known consequence of Kuiper's theorem \cite{Ku} than the (non unital) ring map $\Gamma(X, \B/\A)\to\Gamma(X, M_r(\B/\A))$ induces a bijection between the path components of the associated groups of invertible elements (see for instance \cite{K5} p. 93). Therefore, $\pi_0(\Gamma(X, Fredh (H)))$ may be identified with

$$\lim_{\stackrel{\to}{r}}\pi_0(\Gamma(X, GL_r(\B/\A))) = K_1(\B/\A)$$
and therefore with $K{(\A)}$, as we already mentioned in 2.6.
\end{proof}

\begin{remark}
We may also consider the following ``stabilized'' bundle
$$Fredh_s(H) =\lim_{\stackrel{\to}{n}}Fredh (H^n)$$
and, without Kuiper's theorem, prove in the same way that the set of connected
components of the
 space of sections of this bundle is also isomorphic to $K^{(\A)}(X)$.
\end{remark}

\subsection{} There is an obvious ring homomorphism (since the Hilbert tensor product $H\otimes H$ is isomorphic to $H$)
$$\xymatrix{\K\otimes\K\ar[r]&\K}$$
If $\A$ and $\A'$ are bundles of algebras on $X$ modelled on $\K$, we may use this homomorphism to get a new algebra bundle on $X$, which we denote by $\A\otimes\A'$. From the cocycle point of view, we have a commutative diagram, where the top arrow is induced by complex multiplication

$$\xymatrix{S^1\times S^1\ar[d]\ar[r]&S^1\ar[d]\\
U(H)\times U(H)\ar[r]\ar[d]& U(H\otimes H)\ar[d]\\
PU(H)\times PU(H)\ar[r]& PU(H\otimes H)}$$ It follows that $W_3(\A\otimes \A') =
W_3(\A) + W_3(\A')$ in $H^3(X; \Z)$ and one gets a ``cup-product"

$$\xymatrix{K^{(a)}(X)\times K^{(a')}(X)\ar[r]&K^{(\alpha + {\alpha'})}(X)}$$
(well defined up to non canonical isomorphism: see 2.1). This is a particular case
of a ``graded cup-product" which will be introduced in the next section.

\section{Graded twisted $K$-theory in the finite and infinite-dimensional cases}

\subsection{}

We are going to change our point of view and now consider $\Z/2$-graded
finite-dimensional algebras which are central and simple (in the graded sense). We
are only interested in the complex case. The real case is treated with great details
in \cite{DK} and does not seem to generalize in the infinite-dimensional framework
\footnote{This is not quite true if we work in the context of ``Real'' $K$-theory in
the sense of Atiyah \cite{A1}. We shall not consider this generalization here,
although it looks interesting in the light of
``equivariant twisted $K$-theory'' as we shall show in \S  6.} .\\\\
In the complex case, there are just two ``types" of such graded algebras (up to Morita
equivalence\footnote{This means that we are allowed to take the graded tensor
product with $End(V_0\oplus V_1)$ with the obvious grading.}) which are $\R = \C$ and
$\C\times \C = \C[x]/(x^2 -1)$. For a type $R$ of algebra, the graded inner
automorphisms of $A = R\ootimes End(V_0\oplus V_1)$ may be given by either an
element of degree 0 or an element of degree 1 in $A^*$. This gives us an
augmentation (whose kernel is denoted by $Aut^0(A))$:
$$\xymatrix{Aut (A)\ar[r]&\Z/2}$$
Therefore, for bundles of $\Z/2$-graded algebras modelled on $A$, we already get an
invariant in $H^1(X; \Z/2)$, called the ``orientation" of $A$ and which may be
represented by a real line bundle. A typical example is the (complexified) Clifford
bundle $C(V)$, associated to a real vector bundle $V$ of rank $n$. Its orientation
invariant is the first Stiefel-Whitney class associated to $V$ (cf. \cite{DK}). Note
that the type $R$ of $C(V)$ is $\C$ if $n$ is even and $\C[x]/(x^2-1)\cong\C\times
\C$ if $n$ is odd\footnote{If $V$ is oriented and even dimensional for instance,
this does not imply that $C(V)$ is a bundle of graded algebras of type
$M_2(\Lambda)$ for a certain bundle of ungraded algebras $\Lambda $. However, for
suitable vector bundles $V_0$ and $V_1$, this is the case for the graded tensor
product $C(V)\ootimes End(V_0\oplus V_1)$.}.

\subsection{}

For the second invariant, let us start with $M_{2n}(\C)$ as a basic graded algebra to
fix the ideas and let us put a C*-algebra metric on $A$. We have an exact sequence of
groups as in the ungraded case (where $PU^0(2n) = PU(2n)\cap Aut^0(A))$
$$\xymatrix{
1\ar[r]&
 S^1\ar[r]&
U(n) \times U(n)\ar[r]& PU^0(2n)\ar[r]& 1}$$
Therefore, a bundle with structural
group $PU^0(2n)$ also has a class in $H^2(X; S^1) = H^3(X; \Z)$ which is easily seen
to have order $n$ as in the ungraded case. A similar argument holds if the graded
algebra is $M_{n}(\C)\times M_{n}(\C)$. It follows that the ``graded Brauer group"
$GBr(X)$ is
$$GBr(X) = \Z/2\times H^1(X; \Z/2)\times Tors(H^3(X; \Z))$$
as already quoted in the introduction (see \cite{DK} for the explicit group law on
$GBr(X)$). If $\A$ is a bundle of $\Z/2$-graded finite-dimensional algebras, we
define the graded twisted $K$-theory $K^\A(X)$ as the \underline{graded} $K$-theory
of the graded algebra $\Gamma(X, \A)$ as recalled in Section 1. This definition only
depends on the class of $\A$ in $GBr(X)$. We recover the definition in \cite{DK} by
using again the Serre-Swan theorem as in 2.1. For instance, if we consider the
bundle of (complex) Clifford algebras $C(V)$, associated to a real vector bundle
$V$, the invariants we get are $w_1(V)$ and $W_3(V)$, as quoted in the introduction.
If these invariants are trivial, the bundle $V$ is ${}^c$spinorial of even rank and
$C(V)$ may be identified with the bundle of endomorphisms $End(S^+\oplus S^-)$,
where $S^+$ and $S^-$ are the even and odd ``spinors" associated to the
$Spin^c$-structure.

\subsection{}
In order to define graded twisted $K$-theory in the infinite-dimensional case, we
follow
 the same pattern as in \S  2. For instance, let us take a graded bundle of algebras $\A$
 modelled on $M_2(\K)$: it has 2 invariants, one in $H^1(X; \Z/2)$, the other in
 $H^3(X; \Z)$ (and not just in the torsion part of this group). We then define $K^\A(X)$ as
 the \underline{graded} $K$-theory of the graded algebra $\Gamma(X, \A)$, according to \S  1.
 The same definition holds for bundles of graded algebras modelled on
 $\K\times \K = \K[x]/(x^2 -1) = \K\ootimes\K x$.\\
If $C^{0,1}$ is the Clifford algebra $\C\times \C$ with its usual graded structure,
the general results of \S  1 show that the group $K^\A(X)$ fits into an exact
sequence:

$$\xymatrix@=.6cm{K_1^{(\A\ootimes C^{0,1})} (X)\ar[r]&
K_1^{(\A)} (X)\ar[r]&
K^\A(X)\ar[r]&
K^{(\A\ootimes C^{0,1})} (X)\ar[r]&
K^{(\A)} (X)}$$
where $K_i^{(\A)}(X)$ denotes in general the $K_i$-group of the Banach algebra $\Gamma(X , \A)$.

\subsection{}
Let us now assume that $\A$ is oriented modelled on $M_2(\K)$ (which means that the
structural group of $\A$ may be reduced to $PU^0(H\oplus H)$; see below or 3.2 in
the finite-dimensional situation). We are going to show that $\A$ may be written as
$M_2(\A')$, with the obvious grading, $\A'$ being an ungraded bundle of algebras
modelled on $\K$. For this purpose, we write the commutative diagram
$$\xymatrix{
 S^1\ar[d]\ar[r]&S^1\ar[d]\\
U(H)\ar[r]\ar[d]& U(H)\times U(H)\ar[d]\\
PU(H)\ar[r]&PU^0(H\oplus H)}$$ where the first horizontal map is the identity and
the others are induced by the diagonal. This shows that $H^1(X; PU(H))
\equiv H^1(X; PU^0(H\oplus H))$, which is equivalent to saying that $\A$ may be written as $M_2(\A')$ for a certain bundle of algebras $\A'$.\\
Therefore, $\A\ootimes C^{0,1}$ is Morita equivalent to $\A'\times \A'$ and
$K^\A(X)$ is the $K$-theory of the ring homomorphism (more precisely the functor
defined by the associated extension of the scalars, as we shall consider in other
situations)
$$\xymatrix{\A'\times \A'\ar[r]& M_2(\A')}$$
defined by $(a, b)\mapsto\left(\begin{array}{cc}a&0\\0&b\end{array} \right)$ (no
grading). We should also note that $\A'$ is Morita equivalent to $\A$ as an ungraded
bundle of algebras. Since the usual $K$-theory (resp. graded $K$-theory) is
invariant under Morita equivalence (resp. graded Morita equivalence), the previous
considerations lead to the following theorem:

{\theorem{}{Let $\A$ be an \underline{oriented} bundle of graded algebras modelled on $M_2(\K)$. Then $K^\A(X)$ is isomorphic to $K^{(\A)}(X)$ via the identification above.}}

\subsection{} The same method may be applied in the case when $\A$ is an oriented bundle
of graded algebras modelled on $\K\times \K = \K[x]/(x^2 -1)$. The
\underline{graded oriented} automorphisms of $\K\times \K$ induced by $PU(H\oplus
H)$ are diagonal matrices of type
$$\left(\begin{array}{cc}a&0\\0&a \end{array} \right) $$
This shows that $\A$ is isomorphic to $\A'\times \A'$ and $\A\ootimes C^{0,1}$ is
isomorphic to $M_2(\A')$. Therefore, $K^\A(X)$ is the Grothendieck group of the ring
homomorphism $\A'\to M_2(\A')$, defined by
$$a\mapsto\left(\begin{array}{cc}a&0\\0&a \end{array} \right)$$
Hence we have the following theorem, analogous to 3.5:

{\theorem{}{ Let $\A$ be an \underline{oriented} bundle of graded algebras modelled on $\K\times \K$. Then $\A$ is isomorphic to $\A'\times \A'$ and the group $K^\A(X)$ is isomorphic to $K_1(\A')$ via the identification above.}}

\remark{} One may notice that if $\A$ is a bundle of oriented graded algebras
modelled on $M_2(\K)$, the associated bundle with fiber $M_2(\K^+)$\footnote{$\K^+$ is the ring $\K$ which a unit added} has a section
$\varepsilon$ of degree 0 and of square 1, which commutes (resp. anticommutes) with
the elements of degree 0 (resp. 1); it is simply defined by the matrix
$$\left(\begin{array}{cc}1&0\\0&-1 \end{array} \right)$$

\subsection{} Finally, we would like to give an equivalent definition of $K^\A(X)$ in
terms of Fredholm operators as in \S  2. This was done in \cite{DK} if the class of
$\A$ belongs to the torsion group of $H^3(X; \Z)$ and in \cite{AS1} for the general
case. We shall give another treatment here, using again the $K$-theory of graded
Banach
algebras.\\\\
Following the general notations of \S  2, we have an exact sequence of bundles of
graded Banach algebras
$$\xymatrix{
0\ar[r]&
\A\ar[r]&
\B\ar[r]&
\B /\A\ar[r]&
0}
$$
Since the graded $K$-groups of $\B$ are trivial, we see as in 2.6  that $GrK(\Gamma(X , \A))$ is isomorphic to $GrK_1(\Gamma(X , B /\A))$.\\

In order to shorten the notations, we denote by $B$ the graded Banach algebra
$\Gamma(X , \B)$, by $\Lambda$ the graded Banach algebra $\Gamma(X , B /\A)$ and by
$f$ the surjective map $B\to\Lambda$. The following lemma\footnote{It might be
helpful for a better understanding to notice that the category of finitely generated
free $End(H)$-modules is equivalent to the category of Hilbert spaces $H^n$ for
$n\in\N$. This ``local'' situation is twisted by the group $PU(H)$.} may be proved
in the same way as in \cite{K2} p. 78:

{\lemma{}{Any element of $GrK_1(\Lambda)$ may be written as the homotopy class of a pair
$(E,\varepsilon)$ where $E$ is a free graded $B$-module and $\varepsilon$ is a grading of degree one
of the associated $\Lambda$-module (see 1.4 for the definition of a grading).}}

\subsection{} By the well-known dictionary between modules and bundle theory, we may view
$\varepsilon$ as a grading of a suitable bundle of free $End(H)/\K$-modules. By
spectral theory, we may also assume that $\varepsilon$ is self-adjoint. Finally, following
\cite{K2}, we define a \underline{quasi-grading}\footnote{``quasi-graduation'' in
French.} of $E$ as a family of Fredholm endomorphisms $D$ such that
\begin{enumerate}
\item $D^* = D$
\item $D$ is of degree one
\end{enumerate}
The following theorem is the analogue in the graded case of Theorem 2.8 (cf
\cite{K2} p. 78/79).

{\theorem{}{The (graded) twisted $K$-group $K^\A(X)$ is the Grothendieck group
associated to the semi-group of homotopy classes of pairs $(E, D)$ where $E$ is a
free $\Z/2$-graded $\B$-module and $D$ is a family of Fredholm endomorphisms of $E$
which are self-adjoint and of degree 1.}}

\remark{} Let us assume that $\A$ is oriented modelled on $M_2(\K)$. The description
above gives a Fredholm description of $GrK_1(\A) = K_1(\A )$: we just take the
homotopy classes of sections of the associated bundle of self-adjoint Fredholm
operators $Fredholm^*(\B)$ whose essential spectrum is divided into two non empty
parts\footnote{See the appendix about the role of self-adjoint Fredholm operators in
$K$-theory.}  in ${\R^+}^*$ and ${\R^-}^*$.

\subsection{} This Fredholm description of $K^\A(X)$ enables us to define a cup-product

$$\xymatrix{ K^\A(X)\times K^{\A'}(X)\ar[r]& K^{\A\ootimes \A'}(X)}$$
where $\A\ootimes \A'$ denotes the graded tensor product of $\A$ and $\A'$. This
cup-product is given by the same formula as in \cite{DK} p. 19 and generalizes it:
$$(E, D)\smile (E', D') = (E\ootimes E', D\ootimes 1 + 1\ootimes D')$$
(see the appendix about the origin of this formula in usual topological
$K$-theory).

\subsection{} To conclude this section, let us consider a locally compact space $X$ and
a bundle of graded algebras $\A$ on $X$. For technical reasons, we assume the existence of
a compact space $Z$ containing $X$ as an open subset, such that $\A$ extends to a bundle
(also called $\A$) on $Z$. There is an obvious definition of $K^\A(X)$ as a relative term
in the following exact sequence (where $T = Z -X$ and $\A' = \A\ootimes C^{ 0,1}$):
$$\xymatrix{
K^{\A'}(Z)\ar[r]& K^{\A'}(T )\ar[r]& K^\A(X)\ar[r]& K^\A(Z)\ar[r]& K^\A(T)}$$ By the
usual excision theorem in topological $K$-theory, one may prove that this
definition of $K^\A(X)$ is independent from the choice of $Z$.\\\\
The method described before and also in \cite{K2} \S  3 shows how to generalize the
definition of $K^\A(X)$ in this case: one takes homotopy classes of pairs $(E, D)$
as in 3.12, with the added assumption that the family $D$ is acyclic at infinity. In
other words, there is a compact set $S\subset X$, such that $D_x$ is an isomorphism
when $x\notin S$ (see \cite{K2} p. 89-97 for the technical details of this
approach). This Fredholm description of $K^\A(X)$ will be important in the next
section for the description of the Thom isomorphism.

\section{The Thom isomorphism in twisted $K$-theory\protect\footnote{See also \cite{Ca}.}}

\subsection{} Let $V$ be a finite-dimensional real vector bundle on a locally compact space
$X$ which extends over a compactification of $X$ as was assumed in 3.15. Then the
complexified Clifford bundle $C(V)$ has a well-defined class in the graded Brauer
group of $X$. If $\A$ is another twist on $X$, we can consider the graded tensor
product $\A\ootimes C(V)$ and the associated group $K^{\A\ootimes C(V)}(X)$. As was
described in 3.15, it is the group\footnote{We should note that $E$ is a
$\B$-module, not an $\A$-module. Nevertheless, we shall keep the notation $K^\A$.}
associated to pairs $(E, D)$ where $E$ is a graded bundle of free $\B$-modules and
$D$ is a family of Fredholm endomorphisms which are of degree 1, self-adjoint and
acyclic at $\infty$. Let us now consider the projection $\pi:V\to X$. For
simplicity's sake, we shall often call $E, \A, \B,\dots$ the respective pull-backs
of $E, \A, \B,\dots$ via this projection. Since $\B$ and $C(V)$ are subbundles of
$\B\ootimes C(V)$, $E$ may be provided with the induced $\B$ and $C(V)$-modules
structures.

{\theorem{}{Let $d(E, D) $be an element of $K^{\A\ootimes C(V)}(X)$ with the
notations above. We define an element $t(d(E, D))$ in the group $K^\A(V)$ as
$d(\pi^*(E), D')$ where $D'$ is the family of Fredholm operators on $\pi^*(E)$,
defined over the point $v$ in $V$ (with projection $x$ on $X$) as
$$D'_{(x,v)} = D_x + \rho(v)$$
where $\rho(v)$ denotes the action of the element $v$ of the vector bundle $V$
considered as a subbundle of $C(V)$. The homomorphism
$$t: K^{\A\ootimes C(V)}(X)\to K^\A(V)$$
($t$ for ``Thom") defined above is then an isomorphism.}}

\begin{proof}\footnote{According to a suggestion of J. Rosenberg, it should be
possible to give a proof with the $KK$-theory of Kasparov by describing an explicit
inverse to the homomorphism $t$. However, $KK$-theory is out of the scope of this
paper.} We should first notice that $V$ may be identified with the open unit ball
bundle of the vector bundle $V$ and is therefore an open subset of the closed unit
ball bundle of $V$. Moreover, since $V$ and $\A$ extends to a compactification of
$X$, the required conditions in 3.15 for the
definition of the twisted $K$-theory of $X$ and $V$ are fulfilled.\\\\
We shall now provide two different proofs of the theorem.\\
The first one, more elementary in spirit, consists in using a Mayer-Vietoris
argument which we can apply here since the two sides of the formula above behave as
cohomology theories\footnote{Strictly speaking, one has to ``derive'' the two
members of the formula, which can be done easily since they are Grothendieck groups
of graded Banach categories.} with respect to the base $X$. Therefore, we may assume
that $\A$ and $V$ are trivial: this is a special case of the theorem stated in
\cite{K3} pp. 211/212.\\
The second one is more subtle and may be generalized to the equivariant case. Let us
first describe the Thom isomorphism for complex $V$. This is a slight modification
of Atiyah's argument using the elliptic Dolbeault complex \cite{A2}. More precisely,
we consider the composite map
$$\phi: K^\A(X)\to K^{\A\ootimes C(V)}(X)\stackrel{t}{\to} K^\A(V)$$
The first map $\theta$ is the cup-product with the algebraic ``Thom class" which is
$\Lambda V$ provided with the classical Clifford graded module structure. This map
is an isomorphism from well-known algebraic considerations (Morita equivalence).
Therefore $t$ is an isomorphism if and only if $\phi$ is an isomorphism. On the
other hand, $\phi$ is just the cup-product with the topological Thom class $T_V$
which belongs to the usual topological $K$-theory $K(V)$ of Atiyah and Hirzebruch
\cite{AH2}, \cite{ABS}. In order to prove that $\phi$ is an isomorphism, we may now
use the exact sequence in 3.3 to reduce ourselves to the ungraded twisted case. In
other words, it is enough to show that the cup-product with $T_V$ induces an
isomorphism
$$\Psi:K^{(\A)}(X)\to K^{(\A)}(V)$$
In order to prove this last point, Atiyah defines a reverse map\footnote{More
precisely, one has to define an index map parametrized by a Banach bundle, which is
also classical \cite{FM}.}
$$\Psi':K^{(\A)}(V)\to K^{(\A)}(X)$$
He shows that $\Psi'\Psi=Id$ and, by an ingenious argument, deduces that $\Psi\Psi'= Id$ as
well.\\\\
Now, as soon as Theorem 4.2 is proved for complex $V$, the general case follows from
a trick already used in \cite{K3} p. 241: we consider the following three Thom
homomorphisms which behave ``transitively":

$$\xymatrix@=.6cm{
K^{\A\ootimes C(V)\ootimes C(V)\ootimes C(V)} (X)\ar[r]& K^{\A\ootimes C(V)\ootimes
C(V)} (V)\ar[r]&&\\
&& \!\!\!\!\!\!\!\!\!\!\!\!\!\!\!\!\!\!\!\!\!\!\!\!\!\!\!\!\!\!\!\!\!\!\!\!
K^{\A\ootimes C(V)}(V\oplus V)\ar[r]&K^\A (V\oplus V\oplus V)}$$ We know that the
composites of two consecutive arrows are isomorphisms since $V\oplus V$ carries a
complex structure. It follows that the first arrow is an isomorphism, which is
essentially the theorem stated (using Morita equivalence again).
\end{proof}

\subsection{}
Let $\A$ be a \underline{graded} twist. As we have seen before, it has two
invariants in $H^1(X; \Z/2)$ and in $H^3(X; \Z)$. The first one provides a line
bundle $L$ in such a way that the graded tensor product $\A_1=\A\ootimes C(L)$ is
oriented. From Thom isomorphism and the considerations in 3.5/7, we deduce the
following theorem which gives the relation between the ungraded and graded twisted
$K$-groups.

{\theorem{}{Let $\A$ be a graded twist of type $M_2(\K)$ and $\A_1 = \A\ootimes
C(L)$, where $L$ is the orientation bundle of $\A$. Then $\A_1$ may be written as
$\A_2\times \A_2$ where $\A_2$ is ungraded of type $\K$ . Therefore, we have the
following isomorphisms
$$K^\A(X)\cong
K^{\A\ootimes C(L)\ootimes C(L)}(X)\cong K^{\A_1}(L)\cong K_1^{(\A_2)}(L)$$ Let $\A$
be a graded twist of type $\K\times \K$ . With the same notations, we have the
following isomorphisms
$$K^\A(X)\cong
K^{\A\ootimes C(L)\ootimes C(L)} (X)\cong K^{\A_1}(L)\cong K^{(\A_1)}(L)$$}}

\section{General equivariant $K$-theory}

Note: this section is inserted here for the convenience of the reader as an
introduction to \S  6. It is mainly a summary of results found in \cite{S}, \cite{B},
\cite{LTX} and \cite{AS1}.

\subsection{}
Let $A$ be a Banach algebra and $G$ a compact Lie group acting on $A$ via a
continuous group homomorphism
$$\xymatrix{G\ar[r]&Aut(A)}$$
where $Aut(A)$ is provided with the norm topology. We are interested in the category
whose objects are finitely generated projective $A$-modules $E$ together with a
continuous left action of $G$ on $E$ such that we have the following identity holds
($g\in G$, $a\in A$, $e\in E$)
$$g.(a.e) = (g.a).(g.e)$$
with an obvious definition of the dots. We define $K_G(A)$ as the Grothendieck group of this category $\mathcal{C} = \mathcal{P}_G(A)$. By the well-known dictionary between bundles and modules, we recover the usual equivariant $K$-theory defined by Atiyah and Segal \cite{S} if $A$ is the ring of continuous maps on a compact space $X$ (thanks to the lemma below). It may also be defined as a suitable semi-direct product of $G$ by $A$ (cf. \cite{B}).\\
More generally, if $A$ is a $\Z/2$-graded algebra (where $G$ acts by degree 0
automorphisms), we define the graded equivariant $K$-theory of $A$, denoted by
$GrK_G(A)$, as the Grothendieck group of the forgetful functor
$$\xymatrix{\mathcal{C}^{0,1}\ar[r]&\mathcal{C}}$$
where $\mathcal{C}^{0,1}$ is the category of ``graded objects" in $\mathcal{P}_G(A)$ which is defined as $\mathcal{P}_G(A\ootimes C^{0,1})$.\\
In the same spirit as in \S  1, we define ``derived" groups $K_G^{p,q}(A)$ as
$GrK_G(A\ootimes C^{p,q}$). They satisfy the same formal properties as the usual
groups $K_n(A)$ (also denoted by $K^{-n}(A)$ with $n = q -p$), for instance Bott
periodicity. The following key lemma enables us to translate many general theorems
of $K$-theory into the equivariant framework:

{\lemma{}{Let $E$ be an object of the category $\mathcal{P}_G(A)$. Then $E$ is a
direct summand of an object of type $A\otimes_\C M$ where $M$ is a
finite-dimensional $G$-module.}}

\begin{proof} (\cite{S} p. 134).
Let us consider the union $\Gamma$ of all finite-dimensional invariant subspaces of
the $G$-Banach space $E$. According to a version of the Peter-Weyl theorem quoted in
\cite{S} (loc. cit.), this union is dense in $A$. We now consider a set
$e_1,\dots,e_n$ of generators of $E$ as an $A$-module. Since $\Gamma$ is dense in
$A$ and $E$ is projective, one may choose these generators in the subspace $\Gamma$.
Let $M_1,\dots,M_{n}$ be finite-dimensional invariant subspaces of $E$ containing
$e_1,\dots, e_n$ respectively and let $M$ be the following direct sum
$$M = M_1\oplus\dots\oplus M_{n}$$
We define an equivariant surjection

$$\phi:A\otimes_\C M\cong (A\otimes_\C M_1)\oplus \dots (A\otimes_\C M_{n})\to E$$
between projective left $A$-modules by the formula
$$\phi(\lambda_1\otimes m1,\dots,\lambda_n\otimes m_n) = \lambda_1m_1 +\dots+\lambda_nm_n$$
This surjection admits a section, which we can average out thanks to a Haar measure
in order to make it equivariant. Therefore, $E$ is a direct summand in $A\otimes_\C
M$ as stated in the lemma.\end{proof}

\subsection{}
As for usual $K$-theory, one may also define equivariant $K$-theory for non unital
rings and, using Lemma 5.2 above, prove that any equivariant exact sequence of rings
on which $G$ acts
$$\xymatrix{
0\ar[r]& A'\ar[r]& A\ar[r]& A''\ar[r]& 0 }$$ induces an exact sequence of
equivariant $K$-groups
$$\xymatrix{
K_G^{n-1}(A)\ar[r]&K_G^{n-1}(A'')\ar[r]&K_G^{n}(A')\ar[r]&K_G^{n}(A)\ar[r]&K_G^{n}(A'')
}$$
and an analogous exact sequence in the graded equivariant framework.

\subsection{}
After these generalities, let us assume that $G$ acts on a compact space $X$ and let
$\A$ be a bundle of algebras modelled on $\K$. We define the (ungraded) equivariant
twisted $K$-group $K_G^{(\A)}(X)$ as $K_G(A)$, where $A$ is the Banach algebra of
sections of the bundle $\A$. Similarly, if $\A$ is a bundle of graded algebras
modelled on $\K\times\K$ or $M_2(\K)$, we define the (graded) equivariant twisted
$K$-group $K_G^{\A}(X)$ as $GrK_G(A)$, where $A$ is the graded Banach algebra of
sections of
$\A$.\\
As seen in \S  3, it is also natural to consider free graded $\B$-modules together
with a continuous action of $G$ (compatible with the action on $X$) and a family of
Fredholm operators $D$ which are self-adjoint of degree one, commuting with the
action of $G$. With the same ideas as in \S  3, we can show that the Grothendieck
group of this category is isomorphic to $K_G^{\A}(X)$. This is essentially the
definition proposed in \cite{AS1}.

\subsection{} An example was in fact given at the beginning of the history of twisted $K$-theory. In \cite{DK} \S  8, we defined a ``power operation"
$$P:K^\A(X)\to K_{S_n}^{\A^{\ootimes n}}(X)$$
where $S_n$ denotes the symmetric group on $n$ letters acting on $\A^{\ootimes n}$
(in the spirit of Atiyah's paper on power operations \cite{A1}). According to the
general philosophy of Adams and Atiyah, one can deduce from this $n^{\hbox{th}}$
power operation new ``Adams operations" ($n$ being odd and the group law  in
$GBr(X)$ being written multiplicatively):
$$\Psi^n: K^\alpha(X)\to K^{\alpha^n} (X)\otimes_\Z\Omega_n$$
where $\alpha$ belongs to $\Z/2\times H^1(X; \Z/2)\times H^3(X; \Z)$ and $\Omega_n$
denotes the free $\Z$-module generated by the $n^{\hbox{th}}$ roots of unity in $\C$
(the ring of cyclotomic integers\footnote{The introduction of this ring is
absolutely necessary in the graded case.}). One can prove, following \cite{DK}, that
these additive maps $\Psi^n$ satisfy all the required properties proved by Adams. We
shall come back to them in \S  7, making the link with operations recently defined
by Atiyah and Segal \cite{AS2}.

\subsection{} In order to fix ideas, let us consider the graded version of equivariant twisted $K$-theory $K_G^{\A}(X)$ , where $\A$ is modelled on $\C\times\C$ with the obvious grading. We have again a ``Thom isomorphism":
$$t:K_G^{\A\ootimes C(V)}(X)\to K_G^\A(V)$$
whose proof is the same as in the non equivariant case (using elliptic operators).
On the other hand, if $L$ is the orientation line bundle of $\A$, the group $G$ acts
on $L$ in a way compatible with the action on $X$. Therefore, $\A_1 = \A\ootimes
C(L)$ is an algebra bundle with trivial orientation. and $\A_1\ootimes C^{0,1}$ is
naturally isomorphic to $\A_1\times \A_1 = \A_1[x]/(x^2 -1)$. As in 4.4, this shows
that the graded twisted group $K_G^\A(X)$ is naturally isomorphic to the ungraded
twisted $K$-group $K_G^{\A_1}(L)$. The same method shows that we can reduce graded
twisted $K$-groups to ungraded ones if $\A$ is a bundle of graded algebras modelled
on $M_2(\K)$.

\subsection{} Let us mention finally one of the main contributions of Atiyah and Segal
to the subject (\cite{AS1} \S  6, see also \cite{CKRW}), which we interpret in our
language. One is interested in algebra bundles $\A$ on $X$ modelled on $\K$ ,
provided with a left $G$-action. The isomorphism classes of such bundles are in
bijective correspondence with principal bundles $P$ over $PU(H)$ (acting on the
right) together with a left action of
$G$.\\
Such an algebra bundle $\T$ is called ``trivial" if $\T$ may be written as the bundle of algebras of compact operators in $End(V)$, where $V$ is a $G$-Hilbert bundle. Equivalently, this means that the structure group of $\T$ can be lifted equivariantly to $U(H)$ (in a way compatible with the $G$-action).\\
We now say that two such algebra bundles $\A$ and $\A'$ are equivalent if there exist two trivial algebra bundles $\T$ and $\T''$ such that $\A\otimes\T$ and $A'\otimes\T''$ are isomorphic as $G$-bundles of algebras. The quotient is a group since the dual of $\A$ is its inverse via the tensor product of principal $PU(H)$-bundles. This is the ``equivariant Brauer group" $Br_G(X)$.\\
A closely related definition (in the framework of C*-algebras and for a locally
compact group $G$) is given in \cite{CKRW}. It is very likely that it coincides with
this one for compact Lie groups, in the light of an interesting filtration described
in this paper, probably associated to a spectral sequence.

{\theorem{(cf. \cite{AS1} prop. 6.3)}{ Let $X_G = EG\times_G X$ be the Borel space
associated to $X$. Then the natural map

$$\xymatrix{Br_G(X)\ar[r]&Br(X_G)}$$
is an isomorphism.}}

\subsection{} The interest of this theorem lies in the fact that the equivariant
$K$-theories $K_G(\Gamma(X,\A))$ and $K_G(\Gamma(X, \A'))$ are isomorphic if $\A$
and $\A'$ are equivalent. This follows from the well-known Morita invariance in
operator $K$-theory. We shall study concrete applications of this principle in the
next section.

\section{Some computations of twisted equivariant $K$-groups}

\subsection{} Let us look at the particular case of the ungraded twisted $K$-groups
$K_G^{(\A)}(X)$ where $G$ is a finite group acting on the trivial bundle of algebras
$\A = X\times M_{n}(\C)$ via a group homomorphism $G\to PU(n)$.\footnote{However, we
don't assume that $G$ acts trivially on $X$ in general.} We define $\G$ as the
pull-back diagram
$$\xymatrix{
\G\ar[r]\ar[d]&SU(n)\ar[d]\\
G\ar[r]&PU(n)}$$ Therefore, $\G$ is a central covering of $G$ with fiber $\mu_n$
(whose elements are denoted by Greek letters such as $\lambda$). The following
definition is already present in \cite{K4} \S  2.5 (for $n = 2$):
\begin{definition} A finite-dimensional representation $\rho$ of $\G$ is of ``linear type" if
$\rho(\lambda u) = \lambda\rho(u)$ for any $\lambda\in\mu_n$. We now consider the
category ${E_\G^\A(X)}_l$ whose objects are $\G$-$\A$-modules as before, except that
we request that the $\G$-action be of linear type and \underline{commutes with the
action of $\A$}. By Morita invariance, ${E_\G^\A(X)}_l$ is equivalent to the
category $E_\G(X)_l$ of finite-dimensional $\G$-bundles on $X$, the action of $\G$
on the fibers being of linear type.
\end{definition}

{\theorem{\protect\footnotemark}{The (ungraded) twisted $K$-theory $K_G^{(\A)}(X)$
is canonically isomorphic to the Grothendieck group of the category
$E_{\G}(X)_l$}}\footnotetext{There is an obvious generalization when $\A$ is
infinite-dimensional. However, for our computations, we restrict ourselves to the
finite-dimensional case.}
\begin{proof} One just repeats the argument in the proof of Theorem 2.6 in \cite{K4}, where
$A$ is a Clifford algebra $C(V)$ and $\Z/2$ plays the role of $\mu_n$. We simply
``untwist" the action of $\G$ thanks to the formula (F) written explicitly in the
proof of 2.6 (loc. cit.).\end{proof}

\subsection{}
For $\A = X\times A$ with $A = M_{n}(\C)$, the previous argument shows that
$K_G^{(\A)}(X)$ is a subgroup of the usual equivariant $K$-theory
$K_{\hbox{\scriptsize{$\G$}}}(X)$. From now on, we shall write $K_{G}^{A}(X)$
instead of $K_G^{(\A)}(X)$. Similarly, in the graded case ($A = M_{n}(\C)\times
M_{n}(\C)$ or $M_{2n}(\C)$), we shall write $K_G^A(X)$ instead of $K_G^\A(X)$. If
$X$ is a point and $G$ is finite, $K_G^{(A)}(X)$ is just the $K$-theory of the
semi-direct product $G\ltimes A$.

{\theorem{}{Let $G$ be a finite group acting on the algebra of matrices $A =
M_{n}(\C)$ and let $\G$ be the central extension $G$ by $\mu_n$ described in 6.1.
Then, for $X$ reduced to a point, the group $K_G^{(\A)}(X)= K(G\ltimes A)$ is a free
abelian group of rank the number of conjugacy classes in $G$ which split into $n$
conjugacy classes in $\G$}}

\begin{proof} We can apply the same techniques as the ones detailed in \cite{K4} \S  2.6/12
(for $n = 2$). By the theory of characters on $\G$, one is looking for functions $f$
on $\G$ (which we call of ``linear type") such that
\begin{enumerate}
\item $f(hgh^{-1}) = f(g)$
\item $f(\lambda x) = \lambda f(x)$ if $\lambda$ is an $n^{\hbox{th}}$ root of the unity
\end{enumerate}
The $\C$-vector space of such functions is in bijective correspondence with the
space of functions on the set of conjugacy classes of $G$ which split into $n$
conjugacy classes of $\G$
\end{proof}

\subsection{} Like the Brauer group of a space $X$, one may define in a similar way the
Brauer group $Br(G)$ of a finite group $G$ by considering algebras $A = M_{n}(\C)$
as above with a $G$-action (see \cite{FW} for a broader perspective ; this is also a
special case of the general theory of Atiyah and Segal mentioned at the end of \S  5).
From the diagram written in 6.1, one deduces a cohomology invariant
$$w_2(A)\in H^2(G ; \mu_n)$$
and therefore (via the Bockstein homomorphism) a second invariant $W_3(A)\in H^2(G ;
S^1) = H^2(G ; \Q/\Z) = H^3(G ; \Z)$. It is easy to show that this correspondence
induces a well-defined map
$$W_3:  Br(G)\to H^3(G ; \Z)$$
The following theorem is a special case of 5.9 in a more algebraic situation.

{\theorem{}{ Let $G$ be a finite group. Then the previous homomorphism
$$W_3:  Br(G)\to H^3(G ; \Z)$$
is bijective.}}

\begin{proof} First of all, we remark that $H^3(G ; \Z)\cong H^2(G ; \Q/\Z)$ is the direct
limit of the groups $H^2(G ; \mu_m)$ through the maps $H^2(G ; \mu_m)\to H^2(G ;
\mu_p)$ when $m$ divides $p$. This stabilization process corresponds on the level of
algebras to the tensor product $A\mapsto A\otimes End(V)$, where $V$ is a $G$-vector
space
of dimension $p/m$. Therefore, the map $W_3$ is injective.\\
The proof of the surjectivity is a little bit more delicate (see also 5.9). We can
say first that it is a particular case of a much more general result proved by A.
Fr\"ohlich and C.T.C. Wall \cite{FW} about the equivariant Brauer group of an
arbitrary field $k$: there is a split exact sequence (with their notations)
$$\xymatrix{
0\ar[r]& Br(k)\ar[r]& BM(k, G)\ar[r]& H^2(G ; U(k))\ar[r]& 0 }$$ where $Br(k)$ is
the usual Brauer group of $k$, $U(k)$ is the group of invertible elements in $k$ and
$BM(k, G)$ is a group built out of central simple algebras over $k$ with a
$G$-action. Since $Br(\C) = 0$ and $H^2(G ; U(k)) = H^2(G ; \Q/\Z)$, the theorem is
an immediate consequence.\\

Here is an elementary proof suggested by the referee (for $k = \C$). If we start
with a central extension
$$\xymatrix{
1\ar[r]& \mu_n\ar[r]&\G\ar[r]& G\ar[r]& 1}$$ we consider the finite dimensional
vector space

$$H = \{f \in L^2(\G)\hbox{ such that $f(\lambda\gamma) = \lambda^{-1}f(\gamma)$
with $(\lambda,\gamma)\in\mu_n\times\G$} \}$$ Then $\G$ acts by left translation on
$H$ in such a way that this action is of linear type as described in 6.2. Therefore,
we have a projective representation of $G$ on $H$ and $End(H)$ is a matrix algebra
where $G$ acts.
\end{proof}
\subsection{}
\subsubsection{Remark} Let us consider an arbitrary central extension $\G$
of $G$ by $\mu_n\cong\Z/n$ associated to a cohomology class $c\in H^2(G ; \Z/n)$. We
are interested in the set of elements $g$ of $G$ such that  the conjugacy class
$\langle g\rangle$ splits into $n$ conjugacy classes in $\G_1$. This set only
depends on the image of $c$ in  $H^3(G ; \Z)$ via the Bockstein homomorphism $H^2(G
; \Z/n)\to H^3(G ; \Z)$ (cf. 6.6). In other words, two central extensions of G by
$\Z/n$ with the same associated image by the Bockstein homomorphism have the same
set of $n$-split conjugacy classes.\\\\
In order to show this fact, let us consider the following diagram
$$
\xymatrix{
\mu_n\ar[r]\ar[d]&\mu_{nm}\ar[d]\\
\G_1\ar[r]\ar[d]_\pi&\G\ar[d]\\
G\ar@{=}[r]&G}$$ and an element $g_1$ of $\G_1$. The conjugacy class of $g =
\pi(g_1)$ splits into $n$ conjugacy classes in $\G_1$ if and only if there is a
trace function $f$ on $\G_1$ with values in $\C$ such that $f(g_1c) = f(g_1)c$ when
$c\in\mu_n$. Such a trace function extends obviously to $\G$, which yields to the
result, since the direct limit of the groups $H^2(G ; \Z/nm)$ is precisely $H^2(G ;
\Q/\Z)\cong H^3(G ; \Z)$.

\subsubsection{} This remark may be generalized as follows according to a suggestion
of J.-P. Serre: let $\G_1$ and $\G$ be two group extensions (not necessary central)
of $G$ by abelian groups $C_1$ and $C$ of orders $m_1$ and $m$ respectively, such
that the following diagram commutes (with $\alpha$ injective):
$$\xymatrix{
C_1\ar[r]^\alpha\ar[d]& C\ar[d]\\
\G_1\ar[r]\ar[d]_\pi&\G\ar[d]\\
G\ar@{=}[r]&G}$$ The previous argument shows that if an element $g$ of $G$ splits
into $m$ conjugacy classes in $\G$, it splits into $m_1$ conjugacy classes in
$\G_1$: take trace functions $f$ on $\G$ with values in $C$ such that $f(gc) = f(g)
c$ (we write multiplicatively the abelian group $C$). The converse is true if the
extension of $G$ by $C$ is central.

\subsection{} Let us now assume that $X$ is not reduced to a point. We can use the
Baum-Connes-Kuhn-Slominska Chern character \cite{BC}\cite{Kuh}\cite{SL} which is
defined on $K_\Gamma(X)$ (for any finite group $\Gamma$), with values in the direct
sum $\bigoplus_{\langle\gamma\rangle}H^{\hbox{{\scriptsize even}}}(X^g
)^{C(\gamma)}$. In this formula, $\langle\gamma\rangle$ runs through all the
conjugacy classes of $G$, $C(\gamma)$ being the centralizer of $\gamma$ (the
cohomology is taken with complex coefficients). One of the main features of this
``Chern character"

$$\xymatrix{K_\Gamma(X)\ar[r]&
\bigoplus_{\langle\gamma\rangle}H^{\hbox{{\scriptsize
even}}}(X^\gamma)^{C(\gamma)}}$$ is the isomorphism it induces between
$K_\Gamma(X)\otimes_\Z\C$ and the cohomology with complex coefficients on the
right-hand side. If $E$ is a $\Gamma$-vector
bundle, the map is defined explicitly by  Formula 1.13 p. 170 in \cite{BC}.\\\\
Let us now take for $\Gamma$ the group $\G$ previously considered and let us analyze
the formula in this case. We shall view the right-hand side not just as a function
on the set of conjugacy classes $\langle\gamma\rangle$, but as a function $f$ on the
full group $\Gamma$  with certain extra properties which we shall explain now.\\
If we replace $\gamma$ by $\gamma'$ such that $\gamma' = \sigma\gamma\sigma^{-1}$,
we have a canonical isomorphism $\sigma^*:  H^*(X^{\gamma'})^{C(\gamma')}\to
H^*(X^\gamma )^{C(\gamma)}$ induced by $x\mapsto \sigma x$. This map exchanges
$f(\gamma)$ and $f(\gamma')$ and we have the relation
$f(\gamma)=\sigma^*(f(\gamma'))$,
which says that $f$ is essentially a function on the set of conjugacy classes.\\
On the other hand, if the action of $\G$ is of linear type, we have an extra
relation, an easy consequence of the formula in \cite{BC}, which is $f(\mu\gamma) =
\mu f(\gamma)$ when $\mu$ is an $n^{\hbox{th}}$ root of unity. To summarize, we get
the following theorem.

{\theorem{}{ Let $G$ be a finite group and $A = M_{n}(\C)$ with a $G$-action. Then
the ungraded twisted equivariant $K$-theory $K_G^{(A)}(X)$ is a subgroup of the
equivariant $K$-theory $K_{\G}(X)$, where $\G$ is the pull-back diagram
$$\xymatrix{
\G\ar[r]\ar[d]_\pi&SU(n)\ar[d]\\
G\ar[r]&PU(n)}$$ More precisely, $K_G^{(A)}(X)\otimes_\Z \C$ may be identified with
the $\C$-vector space of functions $f$ on $\Gamma =\G$ with  $f(\gamma)$ in
$H^{\hbox{even}}(X^g)^{C(\gamma)}$ , $\pi(\gamma) = g$, such that the following two
identities hold:
\begin{enumerate}
\item If $\gamma' = \sigma \gamma \sigma^{-1}$, one has $f(\gamma) = \sigma^*(f(\gamma'))$,
according to the formula above.
\item $f(\mu \gamma) = \mu f(\gamma)$ if $\mu$ is an $n^{\hbox{th}}$ root of  unity.
\end{enumerate}
In particular, if $X$  is reduced to a point, we have $\sigma^* = Id$ and
$K_G^{(A)}(X)$ is free with rank the number of conjugacy classes of $G$ which split
into $n$ conjugacy classes in $\G$ (as seen in 6.5.)}}

\subsection{Remarks.} This theorem is not really new. In a closely related context, one
finds similar results in \cite{AR} and \cite{TX}. We should also notice that the
same ideas have been used in \cite{K4} for representations of ``linear type".
Finally, the theorem easily extends to locally compact spaces if we consider
cohomology with compact supports on the right-hand side.

{\theorem{}{  Let $A$ be any finite-dimensional graded semi-simple complex algebra
with a graded action of a finite group $G$. Then the graded $K$-theory
$GrK_0(A')\oplus GrK_1(A')$ of the semi-direct product $A' =G\ltimes A$ is a non
trivial free $\Z$-module. In particular, if $V$ is a real finite-dimensional vector
space with a $G$ action, the group $K_G^A(V)\oplus K_G^A(V\oplus 1)$ is free non
trivial thanks to the Thom isomorphism.}}

\begin{proof} The algebra $G\ltimes A$ is graded semi-simple over the complex numbers.
Therefore, it is a direct sum of graded algebras Morita equivalent to
$M_{n}(\C)\times M_{n}(\C)$ or $M_{2n}(\C)$. In both cases, the graded $K$-theory is
non trivial. The last part of the theorem follows from 4.2.
\end{proof}

The following theorem is a direct consequence of the previous considerations:

{\theorem{}{ Let us now assume that $A = M_{2n}(\C)$ is $G$-oriented as a graded
algebra: in other words, there is an involutive element $\varepsilon$ of $A$ of
degree 0 which commutes with the action of $G$ and commutes (resp. anticommutes)
with the elements of $A$ of degree 0 (resp. 1). Then, the graded $K$-theory
$GrK_*(A')$, with $A' =G\ltimes A$, is a finitely generated free module concentrated
in degree 0. More precisely, one has $GrK_0(A') = K(A')$ and $GrK_1(A') = 0$. In
particular, if $V$ is an even-dimensional real vector space and if $A\ootimes C(V)$
is $G$-oriented, we have (via the Thom isomorphism)

$$K_G^A(V)=K_G^{A\ootimes C(V)}(P)=K(G\ltimes A\ootimes C(V))\hbox{ and }K_G^A(V\oplus 1)=0$$
where $P$ is a point. If we write $A\ootimes C(V)$ as an algebra of matrices
$M_r(\C)$ with a representation $\rho$ of $G$ and call $\G$ the associated central
extension by $\mu_r$, the rank of $K_G^A(V)$ is  the number of conjugacy classes of
$G$ which split into $r$ conjugacy classes in $\G$.}}

In the abelian case, the following two theorems are related to results obtained by
P. Hu and I. Kriz \cite{HK}, using different methods.

{\theorem{}{Let us consider the algebra $A = M_{n}(\C)$ provided with an action of
an abelian group $G$ and $P$ a point. Then the ungraded twisted $K$-theory
$K_G^{(A)}(P)_*= K_*(A')$, with $A' =G\ltimes A$, is concentrated in degree 0 and is
a free $\Z$-module. If we tensor this group with the rationals and if we look at it
as an $R(G)\otimes\Q = \Q[G]$-module, it may be identified  with $R(G')\otimes\Q$
for a suitable subgroup $G'$ of $G$. In particular, the rank of  $K_0(A')$ divides
the order of $G$.}}

\begin{proof} The first part of the theorem is a consequence of the previous more general
considerations. As we have shown before, the algebra $A'$ gives rise to the
following commutative diagram
$$\xymatrix
{\G_n\ar[d]_\pi\ar[r]&SU(n)\ar[d]\\
G\ar[r]& PU(n)}$$ the fibers of the vertical maps being  $\mu_n$. The subset of
elements $\widetilde{g}$ in $\G_n$ such that $\pi(\widetilde{g})$ splits into $n$
conjugacy classes is just the center $Z(\G_n)$ of $\G_n$ (since $G$ is abelian). Let
us put $\Gamma_n = \pi(Z(\G_n))$. Then $K(A')$ may be written as $K_G^A(P)$ where
$P$ is a point. According to Theorem 6.10, this is the subgroup of the
representation ring of $\G_n$ generated by representations of linear type. At this
stage, it is convenient to make $n =\infty$ by extension of the roots of unity, so
that we have an extension of $G$ by $\Q/\Z$
$$\xymatrix{\Q/\Z\ar[r]&\G\ar[r]&G}$$
(the ``linear type" finite-dimensional representations of $\G$ are the same as the
original linear type finite-dimensional representations of $\G_n$). We call
$R(\G)_l)$ the associated Grothendieck group. By the theory of characters, we see
that $R(\G)_l\otimes\Q$ is isomorphic to $R(Z(\G))_l\otimes\Q$, since the characters
of such linear type representations of $\G$ vanish outside $Z(\G)$. On the other
hand,  if we denote by $G'$ the image of $Z(\G)$ in $G$, the extension of abelian
groups
$$\xymatrix{\Q/\Z\ar[r]&Z(\G)\ar[r]&G'}$$
splits (non canonically). This means that we can identify $R(\G)_l\otimes\Q$ with
the representation ring $R(G')\otimes\Q$ as an $R(G)\otimes\Q$-module. This proves
the last part of the theorem.\end{proof}

{\theorem{}{Let us consider a graded algebra $A = M_{2n}(\C)$ provided with a non
oriented action of an abelian group $G$ (with respect to the grading). Then the
graded twisted $K$-theory $K_G^A(P)_*= GrK_*(A')$, with $A' =G\ltimes A$, is
concentrated in a single degree (0 or 1) and is a free $\Z$-module. If we tensor
this group with the complex numbers and if we look at it as an $R(G)\otimes\C =
\C[G]$-module, it may be identified  with $R(G')\otimes\Q$ for a suitable subgroup
$G'$ of $G$. In particular, the rank of $K_G^A(P)_*$ divides the order of $G$.}}
\begin{proof} Let $L$ be the orientation bundle of $A$ (with respect to the action of $G$).
If we change $A$ into $A\ootimes C(L)$ and if we apply the Thom isomorphism theorem,
we have to compute $K_G^A(L)_*$ where $G$ acts on $A$ (resp. $L$) in an oriented way
(resp.
non oriented way).\\
Let us apply Theorem 6.10 in this situation: since $G$ is abelian, the function $f$
of the theorem must be equal to 0 on the elements of $\G$ which are not in $Z(\G)$.
Therefore, the relevant group $K_G^{A}(L)_*$ is reduced to $K_{G'}(L)$ (after
tensoring with $\C$ and where $G'$ is the image of $Z(G)$ in $G$). We now consider
two cases:
\begin{enumerate}
\item the action $\rho$ of $G'$ on $L$ is oriented, in which case we only find $R(G')$ (with
a shift of dimension). Therefore, the dimension of $K_G^{A}(P)\otimes\C$ is the
order of $G'$ which divides the order of $G$.

\item the action $\rho$ of $G'$ on $L$ is not oriented. We then find a direct sum of copies of
$\C$, each one corresponding to an element of $G'$ such that $\rho(g') = -1$. The
dimension of $K_G^A(L)\otimes\C$ is therefore half the order of $G'$, hence divides
$|G|/2$.
\end{enumerate}
\end{proof}

\subsection{Remarks.}
For an oriented action of $G$ on $M_{2n}(\C)$, Theorem 3.5 enables us to solve the
analogous problem in ungraded twisted $K$-theory, which is
done in 6.14.\\
On the other hand, as we already mentioned, the two last theorems are related to
results of P. Hu and I. Kriz \cite{HK}.\\\\
In \cite{K4} we also perform analogous types of computations when $A$ is a Clifford
algebra and $G$ is any finite group, again using the Thom isomorphism. For instance,
if $G$ is the symmetric group $S_n$ acting on the Clifford algebra of $\R^n$ via the
canonical representation of $S_n$ on $\R^n$, there is a nice relation between
$K$-theory and the pentagonal identity of Euler
(cf. \cite{K4} p. 532).\\\\
We would like to point out also that the theory $K_{\pm}(X)$ introduced recently by
Atiyah and Hopkins \cite{AHo} is a particular case of twisted equivariant
$K$-theory. As a matter of fact, it was explicitly present in \cite{K1} \S  3, 40
years ago, before the formal introduction of twisted $K$-theory.
Here is a detailed explanation of this identification.\\
According to \cite{AHo}, the definition of $K_{\pm}(X)$ (in the complex or real
case) is the group $K_{\Z/2}(X\times \mathbb{R}^8)$, where $\Z/2$ acts on $X$ and
also on $\R^8 = \R\times \R^7$ by $(\lambda,\mu)\mapsto(-\lambda, \mu)$. According
to the Thom isomorphism in equivariant $K$-theory (proved in the non spinorial case
in
    \cite{K3}), it coincides with an explicit graded twisted $K$-group $K_{\Z/2}^{A}(X)$,
as defined in \cite{K1}. Here $A$ is the Clifford algebra $C(\R^2)= C^{1,1}$ of
$\R^2$ provided with the quadratic form $x^2 -y^2$ and where $\Z/2$ acts via the
involution $(\lambda,\mu)\to(-\lambda,\mu)$ on $\R\times \R$ (this is also mentioned
briefly in \cite{AHo} p. 2, footnote 1). This identification is valid as well in the
real framework, where we have
8-periodicity.\\
These groups $K_{\Z/2}^{A}(X)$ were considered in \cite{K1} \S  3.3 in a broader
context: $A$ may be any Clifford algebra bundle $C(V)$ (where $V$ is a
\underline{real} vector bundle provided with a non degenerate quadratic form) and
$\Z/2$ may be replaced by any compact Lie group acting in a coherent way on $X$ and
$V$. The paper \cite{K4} gives a method to compute these equivariant twisted
$K$-groups. As quoted in the appendix, the real and complex self-adjoint Fredholm
descriptions (for the non twisted case) which play an important role in \cite{AHo}
were considered independently in \cite{ASi} and \cite{K5}.

\section{Operations on twisted $K$-groups.}

Note:  this section is a partial synthesis of \cite{DK} and \cite{AS2}.

\subsection{} Let us start with the simple case of bundles of  (ungraded) infinite
C*-algebras modelled on $\K$, like in \cite{AS2}. As it was shown in \cite{A1} and
\cite{DK}, we have a $n^{\hbox{{\scriptsize th}}}$ power map
$$P:K^{\A}(X)\to K_{S_n}^{(\A^{\otimes n})}(X)$$
where the symmetric group $S_n$ acts on $\A^{\otimes n}$ by permutation of the
factors.

{\lemma{\protect\footnotemark}{The group $K_{S_n}^{(\A^{\otimes n})}(X)$ is
isomorphic to the group $K_{S_n}^{(\A^{\otimes n})_0}(X)$ where the symbol 0 means
that $S_n$ is acting \underline{trivially} on $\A^{\otimes n}$ }}

\footnotetext{We should note that this lemma is not true for $K_{S_n}(A^{\otimes
n})$ for a general noncommutative ring $A$. Therefore, it is not possible to define
$\lambda$-operations in this case. Twisted $K$-theory is somehow intermediary
between the commutative case and the noncommutative one.}
\begin{proof} As we have shown many times in \S  6, this ``untwisting" of the action of
the symmetric group on $\A^{\otimes n}$ is due to the following fact:  the standard
representation
$$\xymatrix{S_n\ar[r]&PU(H^{\otimes n})}$$
can be lifted into a representation $\rho:  S_n\to U(H^{\otimes n})$ in a way
compatible with the diagonal action of elements of $PU(H)$, a fact which is obvious
to check.\end{proof}

\subsection{Remark.} If we take a bundle of finite dimensional algebras modelled on
$A = End(E)$ where $E = \C^r$, there is another way to check (functorially) this
untwisting of the action of $S_n$ on $A^{\otimes n}$: we identify $A^{\otimes n}$
with $End(E^{\otimes n})$ and $(A^{\otimes n})^*$ with $Aut(E^{\otimes n})$. We have
the following commutative diagram
$$\xymatrix{
&End(E^{\otimes n})^*\ar[d]^\pi&\!\!\!\!\!\!\!\!\!\!\!\!\!\!\!\!\!=Aut(E^{\otimes n})\\
S_n\ar[r]\ar[ru]^\theta&End(E^{\otimes
n})&\!\!\!\!\!\!\!\!\!\!\!\!\!\!\!\!\!=End(E)^{\otimes n} }$$ The vertical map sends
the invertible element $\alpha$ to the automorphism $(u\mapsto\alpha u
\alpha^{-1})$. If $\sigma$ is a permutation, the horizontal map sends $\sigma$ to
the automorphism
$$u_1\otimes\dots\otimes u_n\mapsto u_{\sigma(1)}\dots \otimes u_{\sigma(n)}$$
while the map $\theta$ sends $\sigma$ to the automorphism of $E^{\otimes n}$ defined
by
$$x_1\otimes\dots\otimes x_n\mapsto x_{\sigma(1)}\otimes\dots \otimes  x_{\sigma(n)}$$
Finally, the composition $\pi\theta$, computed on a decomposable tensor of
$E^{\otimes n}$, gives the required result:
$$\xymatrix{
\!\!\!\!\!\!\!\!\!\!\!\!\!\!\!\!\!\!\!\!\!\!\!\!\!\!\!\!\!\! x_1\otimes\dots \otimes
x_n\ar[r]^<(.2)\sigma&
x_{\sigma(1)}\dots \otimes x_{\sigma(n)}\ar[r]^<(.2){u}&
u_1(x_{\sigma(1)})\otimes\dots\otimes u_n(x_{\sigma(n)})
}
$$
$$\xymatrix{
\qquad\qquad\qquad\qquad&\qquad\qquad\qquad\quad& \ar[r]^<(.2){\sigma^{-1}}&
u_{\sigma(1)}(x_1)\otimes \dots\otimes u_{\sigma(n)}(x_n)}$$
\newpage
\subsection{}
\subsubsection{} As was shown in \cite{A1}, a $\Z$-module map
$$\xymatrix{R(S_n)\ar[r]&\Z}$$
defines an operation in twisted $K$-theory by taking the composite of the following
maps
$$\xymatrix@=.6cm{
K^{(\A)}(X)\ar[r]&K_{S_n}^{(\A^{\otimes n})}(X)\cong K_{S_n}^{(\A^{\otimes n})_0}(X)
\cong K^{(\A)^{\otimes n}}(X)\otimes R(S_n) \ar[r]&K^{(\A)^{\otimes n}}(X) }$$ This
is essentially\footnote{The second homomorphism was not explicitly given however.}
what was done in \cite{AS2} \S  10, in order to define the $\lambda^n$ operation of
Grothendieck in this context for instance.

\subsubsection{} Let us call $(F,\nabla)$ or even $F$ (for short) a representative of
the image of $(E , D)$ by the composite of the maps (we now use the Fredholm
description of twisted $K$-theory)
$$\xymatrix{
K^{(\A)}(X)\ar[r]&K_{S_n}^{(\A^{\otimes n})}(X)\ar[r]& K_{A_n}^{(\A^{\otimes
n})_0}(X)}$$ We can define the Adams operations $\Psi^n$ with the method described
in \cite{DK} (which we intend to generalize later on). For this, we restrict the
action of $S_n$ to the cyclic group $\Z/n$ identified with the group of
$n^{\hbox{th}}$ roots of the unity. Let us call $F_r$ the subbundle of $F$ where the
action of $\Z/n$ is given by $\omega^r$ , $\omega$ being a fixed primitive root of
the unity. Then $\Psi^n(E, D)$ is defined by the following sum
$$\Psi^n(E, D)=\sum_{0}^{n-1}F_r\omega^r$$
It belongs formally to $K^{(\A)^{\otimes n}}(X)\otimes_\Z\Omega_n$ where $\Omega_n$
is the ring of $n$-cyclotomic integers. However, if $n$ is prime, using the action
of the symmetric group $S_n$, it is easy to check that $F_r$ is isomorphic to $F_1$
if $r\neq 1$. Therefore we end up in $K^{(\A)^{\otimes n}}(X)$, considered as a
subgroup of $K^{(\A)^{\otimes n}}(X)\underset{\Z}{\otimes}\Omega_n$ , as it was
expected. It is essentially proved in \cite{A1} that this definition of $\Psi^n$
agrees with the classical one.\\\\
There is another operation in twisted $K$-theory which is ``complex conjugation",
classically denoted by $\Psi^{-1}$, which maps $K^{(\A)}(X)$ to $K^{(\A)^{-1}}(X)$
(if we write multiplicatively the group law in $Br(X)$). It is shown in \cite{AS2}
\S 10 how we can combine this operation with the previous ones in order to get
``internal" operations, i.e. mapping $K^{(\A)}(X)$ to itself.

\subsection{}  It is more tricky to define operations in \underline{graded} twisted
$K$-theory. If $\Lambda$ is a $\Z/2$-graded algebra, it is no longer true in general
that a graded involution of $\Lambda$ is induced by an inner automorphism with an
element of degree 0 and of order 2. A typical example is the Clifford algebra
$$(\C\oplus\C)^{\ootimes 2}=(C^{0,1})^{\ootimes 2}$$
which may be identified with the graded algebra $M_2(\C)$. If we put
$$e_1=\left(\begin{array}{cc}0&1\\ 1& 0\end{array}\right)\hbox{ and }
e_2=\left(\begin{array}{cc}0&-i\\i& 0\end{array}\right)$$ we see that there is no
inner automorphism by an element of order 2 and degree 0 permuting
 $e_1$ and $e_2$.

\subsection{}
In order to define such operations in the graded case, we may proceed in at least
two ways.We first remark that if $\A'$ is an oriented bundle modelled on $M_2(\K)$,
the groups $K^{(\A')}(X)$ and $K^{\A'}(X)$ are isomorphic (see 3.5). Moreover, the
ungraded tensor product $\A'^{\otimes n}$ is isomorphic to the graded one
$\A'^{\ootimes n}$ (since $M_2(\C)\ootimes M_2(\C)$ is isomorphic to $M_4(\C)
=M_2(\C)\otimes M_2(\C))$. Let us now assume that the bundle of algebras $\A$ is
modelled on $\K\times\K$ but not necessarily oriented and let $L$ be the orientation
real line bundle of $\A$. Then $\A' = \A \ootimes C(L)$ is oriented modelled on
$M_2(\K)$ ($C(V)$ denotes in general the Clifford bundle associated to $V$). We may
apply the previous method to define operations from $K^{(\A')}(Y) = K^{\A'}(Y)$ to
$K^{\A'^{\ootimes  n}} ( Y)$. If we apply the Thom isomorphism to the vector bundle
$Y = L$ with basis $X$ (see \S  4), we deduce (\underline{for $n$ odd}) operations
from $K^\A(X)$ to $K^{\A^{\ootimes n} }(X)$. However, if $\A$ is oriented, we loose
some information since what we get are essentially the previous operations applied
to the suspension of $X$. Note also that the same method may be applied to
$K^\A(X)$, when $\A$ is modelled on $M_2(\K)$.
\subsection{}
Another way to proceed is to use a variation of the method developed in \cite{DK} p. 21 for
any $\A$ . We consider the following diagram (with $X$ connected)
$$\xymatrix{&F(X, S^1)\ar[r]\ar[d]&F(X, S^1)\ar[d]\\
\Gamma_n\ar[r]^u\ar[d]_\pi&  U^0(\A^{\ootimes n )}
\ar[r]^v\ar[d]&  U^0 (\A\ootimes \A)^{\ootimes n}\ar[d]\\
A_n \ar[r]&  Aut^0( \A^{\ootimes n} ) \ar[r]&  Aut^0 (\A\ootimes \A)^{\ootimes n}}$$
In this diagram $A_n$ is the alternating group (for $n \geq 3$), $\Gamma_n$ is the
cartesian product of $A_n$ and $U^0(\A^{\ootimes  n} )$ over $Aut^0 ( \A^{\ootimes
n} )$ and the horizontal maps between the $U$'s and the $Aut$'s are essentially
given by $g \mapsto  g \otimes  g$. Note that these maps induce a map from $F(X,
S^1)$ to itself
which is $f(z) \mapsto  f(z^2)$.\\
By the general theory developed in 7.2, we know that there exists a canonical map
$\phi$  from $A_n$ to $U^0(\A\ootimes \A)^{\ootimes n}$ such that its composite with
the projection from $U^0 (\A \ootimes \A)^{\ootimes n}$ to $Aut^0  (\A\ootimes
\A)^{\ootimes  n}$ is the composite $\zeta$  of the last horizontal maps. Let us now
define the subgroup $C_n$ of $\Gamma_n$ which elements $x$ are defined by $\varphi
(\pi (x)) = v(u (x))$. It is clear that $C_n$ is a double cover of $A_n$ which is
either the product of $A_n$ by $\Z/2$ or the Schur group $C_n$ \cite{SC} (since
$H^2(A_n ; \Z/2)$ is isomorphic to $\Z/2$ when $n > 3$).\\\\
Using the methods of \S  6, we have therefore the following composite maps
$$\xymatrix{
 K^\A(X)\ar[r]&K_{A_n}^{\A^{\ootimes n}}(X)\ar[r]&
 K^{\A^{\ootimes n}}_{(C_n)}(X) \ar[r]&K^{\A^{\ootimes n}}(X)\otimes
 R(C_n)}$$
where the notation $K_{(C_n )}^{\A^{\ootimes n}} (X)$ means that the group $C_n$
acts trivially on $\A^{\ootimes n}$. Following again Atiyah \cite{A1}, we see then
that any group homomorphism $R(C_n) \to\Z$ gives rise to an operation $K^\A(X)\to
K^{\A^{\ootimes n}} (X)$ Therefore, in the graded case, the Schur group $C_n$
replaces the symmetric group $S_n$, which we have used in the ungraded case.
\subsection{}
In order to define more computable operations, we may replace the alternating group
$A_n$ by the cyclic group $\Z/n$ (when $n$ is odd) as it was already done in
\cite{DK}. The reduction of the central extension $C_n$ becomes trivial and we have
a commutative diagram
$$\xymatrix{
&U^0( \A^{\ootimes n} )\ar[d]\\
\Z/n\ar[ru] \ar[r]&  Aut( \A^{\ootimes n} )
}$$

The Adams operation $\Psi^n$ is then given by the same formula as in 7.4 and
\cite{DK}
$$\Psi^n:  K^\A(X)\to  K^{\A^{\otimes n}} (X)\underset{\Z}{\otimes}  \Omega_n$$
where $\Omega_n$ is the ring of $n$-cyclotomic integers. We can show that this
operation is additive and
multiplicative up to canonical isomorphisms (see Theorem 30, p. 23 in \cite{DK}).\\
We should also notice, following \cite{AS2}, that we can define the Adams operation
$\Psi^{-1}$ in graded twisted $K$-theory as well and combine it with the $\Psi^n$'s
in order to define ``internal'' operations from $K^\A(X)$ to $K^{\A^{\ootimes
n}}(X)\underset\Z\otimes \Omega_n$.

\subsection{} The simplest non-trivial example is

$$\Psi^n: \Z\cong K^1(S^1)\to K^n(S^1)\underset\Z\otimes\Omega_n\cong K^1(S^1)\underset\Z\otimes\Omega_n$$
where $n$ is a product of different odd primes. Since the operation $\Psi^n$ on
$K^2(S^2)$ is the multiplication by $n$, we deduce that
$\theta=\sqrt{(-1)^{(n-1)/2}n}$ belongs to $\Omega_n$ (a well-known result due to
Gauss) and that $\Psi^n$ on $K^1(S^1)$ is essentially the inclusion of $\Z$ in
$\Omega_n$ defined by $1 \mapsto\theta$.

\subsection{}
As a concluding remark, we should notice that the image of $\Psi^n$ as defined in
7.8 is not arbitrary. If $k$ and $n$ are coprimes, the multiplication by $k$ on the
group $\Z/n$ defines an element of the symmetric group $S_n$. The signature of this
permutation is called the Legendre symbol $({}^k_n)$. Moreover, this permutation
conjugates the elements $r$ and $rk$ in the group $\Z/n$. If the Legendre symbol is
1, this permutation can be lifted to the Schur group $C_n$. Let us denote now by
$F_r$ (as in 7.4) the element of the twisted $K$-group associated to the eigenvalue
$e^{2i\pi r}$. Then we see that $F_r$ and $F_{rk}$ are isomorphic if the Legendre
symbol $({}^k_n)$ is equal to 1 since $r$ and $rk$ are conjugate by an element of
the Schur group. If $n$ is prime for instance, $\Psi^n (E)$ may therefore be written
in the following way :
$$\Psi^n (E) = F_0 + \sum_{({}^k_n)=1}U\omega^k+ \sum_{({}^k_n)=-1}V\omega^k$$
where $U$ (resp. $V$) is any $F_k$ with Legendre's symbol equal to 1 (resp. -1).
This shows in particular that the element $\theta = \sqrt{(−1)^{(n−1) /2} n}$ in
7.9 is a ``Gauss sum'', a well-know result.

\section{Appendix: A short historical survey of twisted $K$-theory}
Topological $K$-theory was of course invented by Atiyah and Hirzebruch in 1961
\cite{AH1} after the fundamental work of Grothendieck \cite{BS} and Bott \cite{Bo}.
However, twisted $K$-theory which is an elaboration of it took some time to emerge.
One should quote first the work of Atiyah, Bott and Shapiro \cite{ABS} where
Clifford modules and Clifford bundles (in relation with the Dirac operator) where
used to reinterpret Bott periodicity and Thom isomorphism in the presence of
$\operatorname{Spin}$ or $\operatorname{Spin}^c$ structures. We then started to
understand in \cite{K1}, as quoted in the introduction, that $K$-theory of the Thom
space may be defined almost algebraically as $K$-theory of a functor associated to
Clifford bundles. Finally, it was realized in \cite{DK} that we can go a step
further and consider general graded algebra bundles instead of Clifford bundles
associated to vector bundles (with a suitable metric).\\
On the other hand, the theorem of Atiyah and J\"anich \cite{A1} \cite{J} showed that
Fredholm operators play a crucial role in $K$-theory. It was discovered
independently in \cite{ASi} and \cite{K2}\cite{K5} that Clifford modules may also be
used as a variation for the spectrum of (real and complex) $K$-theory. In
particular, the cup-product from $K^n(X)\times K^p(Y)$ to $K^{n+p}(X\times Y)$ may
be reinterpreted in a much simpler way with this variation. What appeared as a
convenience for usual $K$-theory in these previous references became a necessity for
the definition of the product in the twisted case, as it was shown in
\cite{DK}.\\
The next step was taken by J. Rosenberg \cite{R}, some twenty years later, in
redefining $GBr(X)$ for any class in $H^3(X;\Z)$ (not just torsion classes), thanks
to the new role played by C*-algebras in $K$-theory. The need for such a
generalization was not clear in the 60's (although all the tools were there) since
we were just interested in usual Poincar\'e pairing. In passing, we should mention
that one can define more general twistings in any cohomology theory and therefore in
$K$-theory (see for instance \cite{AS1} p. 18). However, the geometric ones are the
most interesting for the purpose of our paper.\\
We already mentioned the paper of E. Witten \cite{Wi} in the introduction which made
use of this more general twisting of J. Rosenberg, together with \cite{LTX},
\cite{AS1}, \cite{AS2}. But we should also quote many other papers:  \cite{BCMMS}
for the relations with the theory of gerbes and $D$-branes in Physics;  \cite{AR}
for the relation with orbifolds;  \cite{MS} who used the Chern character defined by
A. Connes and M. Karoubi in order to prove an isomorphism between twisted $K$-theory
and a computable ``twisted cohomology'' (at least rationally);  \cite{FHT1} about
the relation between loop groups and twisted $K$-theory;  \cite{FHT2} (again) for
the relation with equivariant cohomology and many other works which are mentioned
inside the paper and in the references. In particular, M. F. Atiyah and M. Hopkins
\cite{AHo} introduced another type of $K$-theory, denoted by them $K_\pm X)$, linked
to deep physical problems, which is a particular case of twisted equivariant
$K$-theory. This definition of $K_\pm (X)$ was already given in \cite{K1} \S  3.3,
in terms of Clifford bundles with a group action (see 6. 16 and also
\cite{K4} for the details).\\
We invite the interested reader to use the classical research tools on the Web for
many other references in mathematics and physics.

\end{document}